\documentclass[12pt]{amsart}
\usepackage{amscd, amssymb}

\newtheorem{pr}{Proposition}
\newtheorem{lm}{Lemma}
\newtheorem{tm}{Theorem}

\newcommand{\proj}{\mathbf P}
\newcommand{\grass}{\mathbf G}
\newcommand{\barr}{\overline}
\newcommand{\rarr}{\rightarrow}
\newcommand{\oh}{{\mathcal{O}}}
\newcommand{\com}{\mathbb{C}}
\newcommand{\Q}{\mathbb{Q}}
\newcommand{\Z}{\mathbb{Z}}
\newcommand{\G}{\mathbf{G}}
\newcommand{\Or}{\mathbf{O}}
\newcommand{\SO}{\mathbf{SO}}
\newcommand{\SL}{\mathbf{SL}}
\newcommand{\GL}{\mathbf{GL}}
\newcommand{\PGL}{\mathbf{PGL}}
\newcommand{\cal}{\mathcal}
\newcommand{\eqq}{\stackrel{\sim}{=}}

\newcommand{\sumo}{\oplus}

\newcommand{\bpf}{\noindent {\em Proof.} }
\newcommand{\epf}{\qed \vspace{+10pt}}

\begin{document}
\begin{center}{\bf The Chow Ring of the Hilbert
Scheme of Rational Normal Curves}
\end{center}
\begin{center}{Rahul Pandharipande$^1$}
\end{center}
\begin{center} {11 June 1996}
\end{center}

\pagestyle{plain}
\footnotetext[1]{Research 
partially supported by an NSF Post-Doctoral Fellowship.}
\setcounter{section}{-1}
\section{\bf{Introduction}}
\subsection{Summary}
Let $\com$ be the ground field of complex numbers.
A rational normal curve in $\proj^d$ is an irreducible, 
nonsingular, non-degenerate, degree $d$ rational curve.
For $d\geq 1$,
let $H(d)$ be the open Hilbert scheme of rational
normal curves of degree $d$ in $\proj^d$.
$H(d)$ is a nonsingular, irreducible, quasi-projective, algebraic
variety. Let $A^*(d)$ be the integral Chow ring of $H(d)$.
In case $d=1$, there is a unique rational normal curve in $\proj^1$.
Hence, $H(1)$ is a point. $H(2)$ is the space
of nonsingular plane conics. 
The dimension of $H(d)$ is $d^2+2d-3$.

In this paper, a presentation
of $A^*(d)$ is computed via the theory of equivariant Chow
groups. The idea is to exhibit $H(d)$ as a quotient of
an appropriate variety $X$ 
by a free $\G$-action. For free actions, the
equivariant Chow ring $A^*_{\G}(X)$ is isomorphic
to the ordinary Chow ring $A^*(d)$ of the quotient $X/\G\eqq H(d)$.
The equivariant Chow ring $A^*_{\G}(X)$ is then computed
in the required cases
via Chow rings of projective bundles and
Chow ideals of degeneracy loci. 

The geometry of
$H(d)$ depends significantly on the parity of $d$.
The quotient approaches and the presentations
of $A^*(d)$ differ for $d$ even and odd. 
In the even case, a $\PGL(2)$-quotient
approach is taken. The geometry of
algebraic $B\PGL(2)$ is studied as a necessary
first step. There is an isomorphism of linear
algebraic groups:
$\PGL(2) \eqq \SO(3)$. The space $B \SO(3)$
is analyzed via conic geometry in projective space.
It is no more difficult to study
algebraic $B \Or(n)$ and $B \SO  (n=2k+1)$ via
higher dimensional quadrics. The equivariant
Chow rings of these two series are computed.
The $d$ odd case
is simpler.
In this case, the quotient
group is taken to be a central extension of
$\SL(2)$ in $\GL(2)$. 
Presentations of $A^*(d)$
in case $d$ is even and odd are determined
in Theorems \ref{heven} and \ref{hodd} respectively.
The equivariant Chow rings of the groups
$\Or(n)$ an $\SO(n=2k+1)$ are computed
in 
Theorem \ref{chor}. 

The equivariant Chow ring
of $\Or(n)$ was first determined by
B. Totaro using 
complex cobordism theory and topology.
After 
the intial algebraic calculation of the ring of $\SO(3)$
presented here, it was realized 
both Totaro's methods and 
the $\SO(3)$ computation generalize to $\SO(2k+1)$.
An algebraic approach to
$B \SO (3)$  is required for the application to
Theorem 1.

In [P],  Chow rings (with $\Q$-coefficients)
of certain moduli spaces of maps are computed via
equivariant Chow groups. The integral
computations presented here were motivated by the
calculations in [P]. These arguments show the
equivariant constructions in [T] and [EG] can
be used effectively to compute ordinary Chow rings
of quotients.

\subsection{Presentations of $A^*(d)$}
\label{prezz}
Equivariant Chow theory is reviewed in section \ref{chow}.
Let $\G$ be a reductive algebraic group. Let $\G\times X \rarr X$
be a linearized algebraic group action on a 
nonsingular quasi-projective
variety $X$. An equivariant Chow ring $A^*_\G(X)$ is defined
via algebraic approximations to $E\G$ and $B \G$.

Let $V$ be a fixed $2$-dimensional $\com$-vector space.
Let $\proj^1\eqq \proj(V)$. There is a canonical isomorphism
$H^0(\proj^1, \oh_{\proj^1}(d)) \eqq Sym^{d}(V^*)$.
Let $$U\subset \bigoplus_{0}^{d} Sym^d(V^*)$$
denote the non-degenerate locus (this is 
the open set consisting of linearly independent  
$(d+1)$-tuples of vectors of $Sym^d(V^*)$).
$U$ parameterizes bases of the linear series of $\oh_{\proj^1}(d)$
on $\proj^1$. 
There is a canonical
$\GL (V)$-action on $U$ with geometric quotient $H(d)$.
The required existence results for the
algebraic quotient problems encountered in this paper
are developed in the Appendix (section \ref{appx}).
$\GL (V)$ acts with finite stabilizers on $U$
(the stabilizer of a point $u\in U$ is the subgroup
of scalar $d^{th}$ roots of unity), By a theorem of
D. Edidin and W. Graham ([EG]), there is a canonical
isomorphism of graded
rings
$$A^*(d) \otimes_{\Z}{\Q} \eqq  A^*_{\GL (V)}(U) \otimes _{\Z}{\Q}.$$ 
The equivariant Chow ring $A^*_{\GL (V)}(U)$ is determined
in section \ref{abe}.
$A^*_{\GL (V)}(U)$ is generated (as a ring) 
in codimensions $1$, $2$ by elements $c_1$, $c_2$ respectively.
There are $d+1$ relations given as follows. Let $S$ be a rank
$2$ bundle with Chern classes $c_1$ and $c_2$. The $d+1$
Chern classes of $Sym^d(S)$ are the relations.
It is not difficult to see $A^i_{\GL (V)}(U) \otimes_{\Z} \Q= 0$
for $i>0$. 
\begin{pr} 
\label{alltor}
$A^*(d)$ is torsion in codimension $i>0$.
\end{pr}
\noindent Note that $\GL(Sym^d (V^*))$ acts transitively
on $U$ and $H(d)$ is a homogeneous space for
$\GL(Sym^d(V^*))=\GL(d+1)$.

Let $\proj(U) \subset \proj(\bigoplus_{0}^{d} Sym^d(V^*))$
be the projective non-degenerate locus. 
$\proj(U)$ is exactly the space of parameterized rational
normal curves.
There is a canonical
$\PGL (V)$-action on $\proj(U)$ with geometric quotient
$H(d)$. This is a free action. Hence, there is a canonical
isomorphism of graded rings (see [EG]):
$$A^*(d)  \eqq  A^*_{\PGL (V)}(\proj(U)).$$ 

Assume $d\geq 2$ is even. Let $d=2n$ (where $n\geq 1$).
$\proj(U) \rarr H(d)$ is a principal $\PGL(V)$-bundle.
Let $S$ be the rank $3$ algebraic vector bundle on $H(d)$
obtained from the principal bundle $\proj(U) \rarr H(d)$
and the representation $Sym^2 (V)$ of $\PGL(V)$.
A discussion of algebraic principal bundles
can be found in the Appendix.
For $1 \leq i \leq 3$, let $c_i\in A^*(d)$ be
the Chern classes of $S$. 
Let $\cal{H}\in A^1(d)$ be the
divisor class of curves meeting a fixed codimension
$2$ linear space in $\proj^d$.
Let $\cal{L}= n \cal{H}$.
In section \ref{evan},
the equivariant Chow ring $A^*_{\PGL (V)}(\proj(U))$ is
evaluated in the even case.
\begin{tm}
\label{heven}
$A^*(d=2n)$ is generated by $c_1$, $c_2$, $c_3$, and $\cal{L}$.
The first relations are:
$$ c_1 =0$$
$$2c_3=0.$$
There are $d+1$ 
additional relations given by the first
$d+1$ Chern classes of the formal expansion:
$$ \frac{(1+\cal{L})^{d+1} \cdot c(Sym^{n-2} (S))}
 {c(Sym^{n} (S))}.$$
(If $n=1$ or $2$, then $c(Sym^{n-2}(S))= 1$.)
\end{tm}
\noindent
It is easily seen from Theorem \ref{heven}
that $A^1(d=2n) \eqq \Z/(d+1)\Z$ with
generator $\cal{L}$. The equation $\cal{L} = n\cal{H}$
can then be uniquely solved to obtain $\cal{H}=2d \cal{L}
=-2 \cal{L}$.

Now assume $d\geq 1$ is odd.
Let $d=2n-1$ (where $n\geq 1$).
Let $$det: \GL (V) \rarr \com^*$$ be the determinant homomorphism.
Let $\Z/n\Z \subset \com^*$ be the subgroup of the $n^{th}$ roots
of unity.
Let $\SL (V,n)= det^{-1}(\Z/n\Z)$.
Consider again $\bigoplus_0^{2n-1} Sym^{2n-1}(V^*)$.
There is a canonical, $\GL (V)$-equivariant, multilinear map
$$\mu:\bigoplus_0^{2n-1} 
Sym^{2n-1}(V^*) \rarr \bigwedge^{2n} Sym^{2n-1}(V^*)$$
given by the exterior product:
$$(\omega_0, \omega_1, \ldots, \omega_{2n-1}) \mapsto
\omega_0 \wedge \omega_1 \wedge \ldots \wedge
\omega_{2n-1}.$$
$\SL (V,n)$ acts trivially on the 1 dimensional space
$\bigwedge^{2n} Sym^{2n-1}(V^*)$. Let 
$Y= \mu^{-1} (p)$ where $0 \neq p \in\bigwedge^{2n} Sym^{2n-1}(V^*)$.
There is an $\SL (V,n)$-action on $Y$. 
In Lemma \ref{freeaq}, it is shown
this is a free action with
geometric quotient $H(d)$.
Hence, there is a canonical
isomorphism of graded rings
$$A^*(d=2n-1)  \eqq  A^*_{\SL (V,n)}(Y).$$ 

Let $S$ now denote the rank $2$ algebraic
vector bundle obtained from the principal
$\SL(V,n)$-bundle $Y \rarr H(d)$ and the
standard representation $V$.
For $1\leq i \leq 2$, let $c_i\in A^*(d)$
be the Chern classes of $S$.
The equivariant Chow ring $A^*_{\SL (V,n)}(Y)$ is evaluated
in section \ref{ode}. 

\begin{tm}
$A^*(d=2n-1)$ is generated by $c_1$ and $c_2$.
The first relation is
$$ nc_1 =0.$$
There are $d+1$ 
additional relations given by the first
$d+1$ Chern classes of $Sym^d(S)$.
\label{hodd}
\end{tm} 

\noindent
It is easily seen that $A^1(d=2n-1) \eqq \Z/n\Z$.

\subsection{Chow rings of the Orthogonal Groups}
The Chow ring of a reductive algebraic group $\G$ is, by
definition, the equivariant Chow ring $A^*_\G(
\text{point})$.
Let $\Or (n)$ and $\SO (n)$ denote the orthogonal and
special orthogonal algebraic groups.
The equivariant calculations of Theorem \ref{heven}
require knowledge of $B\PGL (2)$. $\PGL (2)$ is isomorphic
to $\SO (3)$. The following Theorem will be established:
\begin{tm}
\label{chor}
The integral Chow ring of $\Or (n)$ is    
generated by the Chern classes $c_1, \ldots, c_n$
of the standard representation. The odd classes
are 2-torsion:
$$A^*_{\Or (n)}(\text{point})= 
\Z[c_1, \ldots, c_n]/(2c_1, 2c_3, 2c_5, \ldots).$$

\noindent The integral Chow ring of $\SO (n=2k+1)$
is generated by the Chern classes $c_1, \ldots, c_n$
of the standard representation. The odd classes
are 2-torsion and $c_1=0$:
$$A^*_{\SO (n)}(\text{point})=
 \Z[c_1, \ldots, c_n]/(c_1, 2c_3, 2c_5, \ldots).$$
\end{tm}
\noindent
The Chow ring of $\SO (2k)$ is not generated by the
Chern classes of the standard representation. The
main difference in the odd and even cases is that
$\SO (2k+1) \eqq \Z/2\Z \times \Or(2k+1)$ while
such a product decomposition does not hold
for $\SO (2k)$.
The methods
of this paper do not yield a computation of $A^*_{\SO (2k)}(
\text{point})$. The Chow ring of 
$\SO(n)$ has been computed with $\Q$-coefficients in [EG2].

\subsection{Acknowledgments} 
Equivariant Chow groups were first defined in [T].
Thanks are due to D. Edidin, W. Graham, and B. Totaro
for conversations in which the theory of
equivariant Chow groups was explained. 
The author particularly wishes to thank B. Totaro
for his insights on $\Or(n)$ and $\SO(n)$.
Discussions with
W. Fulton on many related issues have also been helpful.

\section{\bf{Chow Ideals of Degeneracy Loci}}
\subsection{Presentations}
\label{idealz}
For the Chow computations in this paper,
presentations of four ideals associated to tautological 
degeneracy loci are needed. 

Let $E$ be a rank $e$ vector bundle on a 
nonsingular algebraic variety $M$.
We will consider two affine and
two projective fibrations over $M$:
\begin{enumerate}
\item[(i)] $\sumo_{1}^e E \rarr M$,
\item[(ii)] $\proj(\sumo_{1}^e E) \rarr M$,
\item[(iii)] $Sym^2 E^* \rarr M$,
\item[(iv)] $\proj(Sym^2 E^*) \rarr M$.
\end{enumerate}
The 
subspace projectivization is taken in (ii) and (iv).
Let $r=e^2$ denote the
rank of $\sumo_{1}^{e} E$. 
Let $L$ in $A^1(\proj(\sumo_{1}^{e} E))$ be the
class of $\oh_{\proj}(1)$ obtained from the
projectivization. The Chow ring of $\proj(\sumo_{1}^{e} E)$
has a standard presentation:
$$A^*(M)[L]\ /\ \big(L^{r}+
 c_1(\sumo_{1}^{e} E)\cdot 
L^{r-1} + \ldots +c_{r}(\sumo_{1}^{e} E)\big).$$
Similarly, the Chow ring of $\proj(Sym^2 E^*)$ 
has a presentation:
$$A^*(M)[L]\ /\ \big(L^{s}+
 c_1(Sym^2 E^*)\cdot 
L^{s-1} + \ldots +c_{s}(Sym^2 E^*)\big)$$
where $s=\frac{1}{2}(e^2+e)$ is the rank of
$Sym^2 E^*$ and $L$ is again  the class of $\oh_{\proj}(1)$
obtained from the projectivization.
The Chow rings of the affine fibrations (i) and
(iii) are canonically isomorphic to $A^*(M)$.

There are intrinsic, fiberwise degeneracy loci
in these fibrations.
Let  $D_1 \subset \sumo_{1}^{e} E$ and 
$\proj(D_1) \subset \proj(\sumo_{1}^{e} E)$ 
be the
closed subvariety of linearly dependent $e$-tuples
of vectors in the fibers of $E$.
Let $D_2 \subset Sym^2 E^*$ and
$\proj(D_2) \subset \proj(Sym^2 E^*)$
be the closed subvariety of degenerate quadratic forms
on the fibers of $E$.
Let $$I_1\subset A^*(\sumo_{1}^{e} E)\eqq A^*(M), \ \
 J_1\subset A^*(\proj(\sumo_{1}^{e} E)),$$
$$I_2 \subset A^*(Sym^2 E^*)\eqq A^*(M), \ \ 
J_2 \subset A^*(\proj(Sym^2 E^*))$$ be the
ideals generated by classes supported on the
degeneracy loci $D_1$, $\proj(D_1)$, $D_2$,
and $\proj(D_2)$ respectively.
In this section, simple sets of generators of the
ideals $I_1$, $J_1$, $I_2$, and $J_2$ are determined. 

The results of this section are essentially
 special cases of 
Pragacz's presentations of the ideals of Chow classes supported
on degeneracy loci of bundle maps
([Pr]). Pragacz considers more general degeneracy
loci and obtains presentations of their
universal Chow ideals via Schur $S$-polynomials.
Actual (not universal) Chow ideal presentations
are needed here. Since the geometry of the cases
(i)-(iv) is particularly simple, the actual and
the universal presentations coincide.
A full proof will be given here.

For a rank $f$ bundle $F$, let
$c(F)=1+c_1(F)+ \ldots +c_f(F)$.
\begin{lm}
\label{petey}
$I_1\subset A^*(M)$ is generated by
$(\alpha_1, \ldots, \alpha_e)$ where
$$\frac{1}{c(E^*)}=
1+ \alpha_1 +  \ldots + \alpha_e+ \ldots.$$
\end{lm}
\begin{lm}
\label{pete}
$J_1\subset A^*(\proj(\sumo_{1}^{e} E))$ is generated by
$(\alpha'_1, \ldots, \alpha'_e)$ where
$$\frac{c(\sumo_{1}^{e} \oh_{\proj}(1))}{c(E^*)}=
1+ \alpha'_1 +  \ldots + \alpha'_e+ \ldots.$$
\end{lm} 
\begin{lm}
\label{paul}
$I_2\subset A^*(M)$ is generated by
$(\beta_1, \ldots, \beta_e)$ where
$$\frac{c(E^*)}{c(E)}=
1+ \beta_1 +  \ldots + \beta_e+ \ldots.$$
\end{lm}
\begin{lm}
\label{pauly}
$J_2\subset A^*(\proj(Sym^2 E^*))$ is generated by
$(\beta'_1, \ldots, \beta'_e)$ where
$$\frac{c(E^*\otimes \oh_{\proj}(1))}{c(E)}=
1+ \beta'_1 +  \ldots + \beta'_e+ \ldots.$$
\end{lm}
\noindent
The proofs of
Lemmas \ref{petey} -- \ref{pauly} 
are essentially the same. 
The first step is to find a tower of bundles dominating
the degeneracy loci $D_1$, $\proj(D_1)$, $D_2$,
and $\proj(D_2)$.

First consider $D_1$ and $\proj(D_1)$. 
Let $\eta:\proj(E^*) \rarr M$ be the projective bundle.
A point $\xi \in \proj(E^*)$ is a pair $(m,h)$ where
$m\in M$ and $h \in \proj(E^*_m)$.
Let $B$ be the vector bundle on $\proj(E^*)$ determined
as follows. The fiber of $B$ at the point $(m,h)$
is the linear subspace of $\sumo_{1}^{e} E_m$
consisting of $e$-tuples of vectors annihilated by
$h$.
$B$ is a sub-bundle of  $\eta^*(\sumo_{1}^{e} E)$.
There are canonical, proper, surjective projections:
$$\rho: B \rarr D_1 \subset \sumo_{1}^{e} E,$$
$$\proj(\rho): \proj(B) \rarr \proj(D_1)
\subset \proj(\sumo_{1}^{e} E).$$
There are  stratifications of $D_1$ and $\proj(D_1)$
by the
rank of the span of the $e$-tuple of vectors. Over these strata,
$\rho$ and $\proj(\rho)$  are projective bundles.
Hence $\rho$ and $\proj(\rho)$ induce {\em surjections}
on the integral Chow rings
via push-forward:
$$\rho_*: A^*(B) \rarr A^*(D_1),$$
$$\proj(\rho)_*: A^*(\proj(B)) \rarr A^*(\proj(D_1)).$$

Lemmas \ref{petey} and \ref{pete}
are proven by computing the images of
the generators of $A^*(B)$ and $A^*(\proj(B))$ 
respectively. 
Consider the commuting diagrams:
\begin{equation}
\label{heyhey}
\begin{CD}
 B @>{\rho}>>  D_1 \subset \sumo_{1}^{e} E \\
@V{\pi}VV @VVV \\
\proj(E^*) @>{\eta}>> M \\
\end{CD}
\end{equation}
\begin{equation}
\label{heyheyy}
\begin{CD}
 \proj(B) @>{\proj(\rho)}>> \proj( D_1)
 \subset \proj(\sumo_{1}^{e} E) \\
@V{\pi}VV @VVV \\
\proj(E^*) @>{\eta}>> M \\
\end{CD}
\end{equation}
$A^*(B)$ is generated over $A^*(M)$ by 
the class corresponding to $\oh_{\proj(E^*)}(1)$.
 Let
this class be denoted by $\zeta$.
It follows that
$$I_1=(\rho_{*}(1), \rho_{*}(\zeta^1), \rho_{*}(\zeta^2),
\ldots, \rho_{*}(\zeta^{e-1})).$$
Similarly, there is a presentation of $J_1$:
$$J_1=(\proj(\rho)_*(1), \proj(\rho)_{*}(\zeta^1), 
\proj(\rho)_{*}(\zeta^2),
\ldots, \proj(\rho)_{*}(\zeta^{e-1})).$$
To prove Lemmas \ref{petey} and \ref{pete}, 
it is sufficient to establish
the equalities
\begin{equation}
\label{eqql}
 \rho_{*}(\zeta^{i-1})= \alpha_i, \ \ 
\proj(\rho)_{*}(\zeta^{i-1})= \alpha'_i
\end{equation}
for $1 \leq i \leq e$.

First the equalities (\ref{eqql}) for Lemma
\ref{petey} are proven. By definition,
$B\subset \eta^*(\sumo_{1}^{e} E)$. 
In fact, there is a natural exact squence on $\proj(E^*)$:
\begin{equation}
\label{yess1}
0 \rarr B \rarr \eta^*(\sumo_{1}^{e} E) \rarr \sumo_{1}^{e}
\oh_{\proj(E^*)}(1) \rarr 0.
\end{equation}
As a first
step, the class of $[B] \in A^*(\eta^*(\sumo_{1}^{e} E))$
is computed. Since
$\eta^*(\sumo_{1}^{e} E)$
is a projective bundle over $\sumo_{1}^{e} E$, 
$A^*(\eta^*(\sumo_{1}^{e} E))$
is generated over $A^*(M)$ by $\zeta$ (which
satisfies the Chern relation).
By sequence (\ref{yess1}) and Lemma \ref{fullt} below,
it follows that $[B] = \zeta ^{e} \in A^*(\eta^*(\sumo_{1}^{e} E))$.
Denote the natural projection $\eta^*(\sumo_{1}^{e} E)) 
\rarr \sumo_{1}^{e}
E$ by $\phi$.
There is a fundamental equality:
$$\rho_{*}(\zeta ^{i-1}) = \phi_{*} (\zeta^{i-1} \cap [B]) \in A^*(M).$$
The right side is easy to calculate.
$$\phi_{*}(\zeta^{i-1}\cap [B])= \phi_{*}(\zeta^{e-1+i}).$$
For $1\leq i \leq e$, the latter is simply the $i^{th}$ Segre
class of $E^*$. Lemma \ref{petey} is proved.

Lemma \ref{pete} is only slightly more complicated.
The class of $[\proj(B)] \in A^*(\proj(\eta^*(\sumo_{1}^{e} E)))$
is computed. Again 
$A^*(\proj(\eta^*(\sumo_{1}^{e} E)))$
is generated over $A^*(\proj(\sumo_{1}^{e} E))$ by $\zeta$. 
By sequence (\ref{yess1}) and Lemma \ref{fullt},
it follows that $[\proj(B)] 
= (L+\zeta) ^{e} \in A^*(\proj(\eta^*(\sumo_{1}^{e} E)))$
where $L$ is the class of $\oh_{\proj(\sumo_{1}^{r} E)}(1)$.
Denote the natural projection $\proj(\eta^*(\sumo_{1}^{e} E))
 \rarr \proj(\sumo_{1}^{e}
E)$ by $\proj(\phi)$.
There are equalities:
$$\proj(\rho)_{*}(\zeta ^{i-1}) 
= \proj(\phi)_{*} (\zeta^{i-1} \cap [B]) 
=\proj(\phi)_{*}
(\zeta^{i-1}\cap (L+\zeta)^e)$$
in $A^*(\proj
(\sumo_{1}^{e} E))$.
Lemma \ref{ssegre} now yields Lemma \ref{pete}.

The degeneracy loci
$D_2$ and $\proj(D_2)$ are considered next. The notation will
parallel the notation used in the proofs of Lemmas \ref{petey}
and \ref{pete}.
Let $\eta:\proj(E) \rarr M$ be the projective bundle.
A point $\xi \in \proj(E)$ is a pair $(m,p)$ where
$m\in M$ and $p \in \proj(E_m)$.
Let $B$ be the vector bundle on $\proj(E)$ determined
as follows. The fiber of $B$ at the point $(m,p)$
is the linear subspace of quadratic forms on $E_m$
singular at $p$.
$B$ is a sub-bundle of  $\eta^*(Sym^2 E^*)$.
There are  canonical, proper, surjective projections:
$$\rho: B \rarr D_2 \subset Sym^2 E^*,$$
$$\proj(\rho): \proj(B) \rarr \proj(D_2) \subset Sym^2 E^*.$$
There are  stratifications of $D_2$ and $\proj(D_2)$
by the
rank of the quadratic form. Over these strata,
$\rho$ and $\proj(\rho)$  are projective bundles.
Hence $\rho$ and $\proj(\rho)$ induce {\em surjections}
on the integral Chow rings
via push-forward:
$$\rho_*: A^*(B) \rarr A^*(D_2),$$
$$\proj(\rho)_*: A^*(\proj(B)) \rarr A^*(\proj(D_2)).$$

Lemmas \ref{paul} and \ref{pauly}
are proven by computing the images of
the generators of $A^*(B)$ and $A^*(\proj(B))$ 
respectively. 
As before,
$$I_2=(\rho_{*}(1), \rho_{*}(\zeta^1), \rho_{*}(\zeta^2),
\ldots, \rho_{*}(\zeta^{e-1})),$$
$$J_2= (\proj(\rho)_{*}(1), \proj(\rho)_{*}
(\zeta^1), \proj(\rho)_{*}(\zeta^2),
\ldots, \proj(\rho)_{*}(\zeta^{e-1})).$$
To prove Lemma \ref{paul} and \ref{pauly}, 
it is sufficient to establish
the equalities
\begin{equation}
\label{eqql2}
 \rho_{*}(\zeta^{i-1})= \beta_i, \ \ 
\proj(\rho)_{*}(\zeta^{i-1})= \beta'_i
\end{equation}
for $1 \leq i \leq e$.

First the equalities (\ref{eqql2}) for Lemma
\ref{paul} are proven. There 
is an exact squence on $\proj(E)$:
\begin{equation}
\label{yess2}
0 \rarr B \rarr \eta^*(Sym^2 E^*) \rarr E^*\otimes \oh_{\proj(E)}(1)
 \rarr 0.
\end{equation}
The class of $[B] \in A^*(\eta^*(\sumo_{1}^{e} E))$
is computed.  
$A^*(\eta^*(Sym^2 E^*))$
is generated over $A^*(M)$ by $\zeta$.
By sequence (\ref{yess2}) and Lemma \ref{fullt} below,
it follows that 
$$[B] = c_e (E^*\otimes \oh_{\proj(E)}(1) )\in A^*(\eta^*(Sym^2 E)).$$
Denote the natural projection $\eta^*(Sym^2 E^*)) 
\rarr Sym^2
E^*$ by $\phi$.
There is an equality:
$$\rho_{*}(\zeta ^{i-1}) = \phi_{*} (\zeta^{i-1} \cap [B]) \in A^*(M).$$
Lemma \ref{ssegre} now  yields Lemma \ref{paul}.

Lemma \ref{pauly} is established next.
By sequence (\ref{yess2}) and Lemma \ref{fullt} below,
it follows that 
$$[\proj(B)] = c_e \bigg(
\frac{c(E^*\otimes \oh_{\proj(E)}(1))}
{c(\oh_{\proj(Sym^2 E^*)}(-1))}  \bigg)\in 
A^*(\proj(\eta^*(Sym^2 E^*))).$$
There is an equality (since $E^*\otimes \oh_{\proj(E)}(1)$
is a rank $e$ bundle):
$$c_e \bigg(
\frac{c(E^*\otimes \oh_{\proj(E)}(1))}
{c(\oh_{\proj(Sym^2 E^*)}(-1))}  \bigg)=
c_e(E^*\otimes \oh_{\proj(E)}(1) \otimes \oh_{\proj(Sym^2 E^*)}(1)).$$
Denote the natural projection $\proj(\eta^*(Sym^2 E^*)) 
\rarr \proj(Sym^2
E^*)$ by $\proj (\phi)$.
There is an equality:
$$\proj(\rho)_{*}(\zeta ^{i-1}) = 
\proj(\phi)_{*} (\zeta^{i-1} \cap [B]) \in A^*(\proj(Sym^2 E^*)).$$
Lemma \ref{ssegre} now  yields Lemma \ref{pauly}.

\subsection{Lemmas}
The following Lemmas were used in the proofs
of Lemmas \ref{petey} -- \ref{pauly}.
Let $F \rarr N$ be a vector 
bundle on a nonsingular algebraic
variety $N$. 
\begin{lm}
\label{fullt}
Let
$\ 0 \rarr B \rarr F \rarr Q \rarr 0$
be an exact sequence of bundles on $N$.
Let $q$ be the rank of $Q$.
The class $[B]\in A^*(F)\eqq A^*(N)$ is determined
by
$$[B]= c_q(Q).$$
The class $[\proj(B)] \in A^*(\proj(F))$ is determined
by
$$[\proj(B)]= c_q \bigg(\frac {c(Q)}{c(\oh_{\proj(F)}(-1))}
\bigg).$$
\end{lm}
\bpf
This is an  application of the Thom-Porteous formulas 
for degeneracy loci of bundle maps 
(see [F]).
\epf

\noindent
Let $f$ be the rank $F$.
Let $\phi: \proj(F) \rarr N$ be the projection.
\begin{lm}
\label{ssegre}
Let $G$ be a bundle of rank $g=f$ on $N$.
Let $\zeta= c_1(\oh_{\proj(F)}(1))$.
Let $\gamma_{i}$ 
be determined by
$$\frac{c(G)}{c(F)} = 1 + \gamma_{1} + \ldots+
\gamma_{f} + \ldots.$$
Then, for $1 \leq i \leq f$,
 $\gamma_{i}=  \phi_{*}\big( \zeta^{i-1} \cap 
c_f(G \otimes
\oh_{\proj(F)}(1))\big).$
\end{lm} 
\bpf
A simple Segre class argument yields the result.
\epf

\section{\bf{Equivariant Chow Groups}}
\label{chow}
Let $\G$ be a group. 
Let $\G\times X \rarr X$ be a left group action.
In topology, the $\G$-equivariant cohomology of 
$X$ is defined as follows. Let $E\G$ be a contractible
topological space equipped with a free left $\G$-action and 
quotient $E\G/\G=B\G$.
Consider the left action of $\G$ on $X\times E\G$ defined
by:
$$g(x,b)= (g(x), g(b)).$$
$\G$ acts freely on $X\times E\G$. Let
$X\times^{\G} E\G$ be the (topological) quotient. 
The $\G$-equivariant cohomology of 
of $X$, $H_\G^*(X)$, is defined by:
$$H_\G^*(X) = H^*_{sing}(X\times^{\G} E\G).$$
If $X$ is a  locally trivial principal $\G$-bundle,
then $X\times^{\G} E\G$ is a locally trivial
fibration of $E\G$ over the quotient $X/\G$.
In this case, $X\times^{\G} E\G$ is homotopy equivalent
to $X/\G$ and
$$H_\G^*(X) = H^*_{sing}(X\times^{\G} E\G) \eqq 
H^*_{sing}(X/\G).$$
For principal bundles, computing the
equivariant cohomology ring is equivalent
to computing the cohomology of the quotient.

There is an analogous equivariant theory of Chow groups
developed by B. Totaro in case $X$ is a point and
generalized by D. Edidin and W. Graham to arbitrary $X$
([T], [EG]). Let $\G$ be a reductive algebraic group.
Let $\G\times X \rarr X$ be a linearized  algebraic $\G$-action.
The algebraic analogue of $E\G$ is attained by
approximation. Let $V$ be a $\com$-vector space. Let
$\G\times V \rarr V$ be an
algebraic representation of $\G$. 
Let $W\subset V$ be a $\G$-invariant open set satisfying:
\begin{enumerate}
\item[(i)] The complement of $W$ in $V$ is of codimension greater than 
$q$.
\item[(ii)] $\G$ acts freely on $W$ (see the Appendix for
the definition).
\item[(iii)] There exists a geometric quotient $W\rarr W/\G$.
\end{enumerate}
$W$ is an approximation of $E\G$ up to codimension $q$.
By (iii) and the
assumption of linearization, 
a geometric quotient $X\times ^{\G} W$ exists as
an algebraic variety. Let $d=dim(X)$, $e=dim( X\times ^{\G} W)$.
The equivariant Chow groups are defined by:
\begin{equation}
\label{defff} 
A^{\G}_{d-j}(X)= A_{e-j}(X\times ^{\G} W)
\end{equation}
for $0\leq j \leq q.$
An argument is required to check these equivariant
Chow groups are well-defined (see [EG]).
The basic functorial properties of equivariant
Chow groups are  
established in [EG]. In particular, if $X$ is
nonsingular, there is a natural intersection
ring structure on $A_i^{\G}(X)$.

Let $Z$ be a variety of dimension $z$.
For notational convenience, a superscript will
denote the Chow group codimension:
 $$A^{\G}_{z-j}(Z) = A^j_\G(Z), \ A_{z-j}(Z)=A^j(Z).$$
In particular, equation (\ref{defff}) becomes:
$$\forall\  0\leq j \leq q, \ \
A_{\G}^{j}(X)= A^j(X\times ^{\G} W).$$
The following result of [EG] 
will be used.
\begin{pr}
\label{dane}
Let $\com$ be the ground field of complex numbers.
Let $X$ be a quasi-projective variety. Let $\G$ be a reductive
group.
Let $\G\times X \rarr X$ be a linearized proper $\G$-action.
Let $X\rarr X/\G$ be a
quasi-projective
geometric quotient. 
\begin{enumerate}
\item[(i)] If the action is free, then there is a canonical
isomorphism of graded rings:
$$A^*_\G(X) \eqq A^*(X/\G).$$
\item[(ii)] If $\G$ acts with finite stabilizers on $X$, then
there is a canonical isomorphism of graded rings:
$$ A^*_\G(X) \otimes \Q \eqq A^*(X/\G) \otimes \Q.$$
\end{enumerate}
\end{pr}
\noindent Proposition \ref{dane} is a characteristic 0
specialization of Theorem 2 of [EG].

\section{\bf The Chow Rings of $\Or (k)$ and $\SO (2k+1)$}
\subsection{$B\Or (V)$ and $B\SO (V)$}
\label{orthoo}
Let $V$ be a complex vector space equipped with
a non-degenerate quadratic form. Let $\Or (V)$, $\SO (V)$
be the orthogonal and special orthogonal groups
respectively.
Approximations to $E\Or (V)$ and $E\SO (V)$ are
obtained via direct sums of the 
representation $V^*$.
Let $m>>0$ and let
$$W_m \subset \sumo_{1}^{m} V^*$$
denote the spanning locus. $W_m$ is the
locus of $m$-tuples of vectors of $ V^*$
which span $V^*$.
The natural actions of $\Or (V)$ and $\SO (V)$ on
$W_m$ are free and have a geometric quotients
(see section \ref{appx}).
The codimension of the complement of
$W_m$ in $\sumo_{1}^{m}  V^*$ is $m-dim(V^*)+1$.
$W_m$ is an approximation of $E\Or (V)$ and
$E\SO (V)$ up to codimension $m-dim(V^*)$.
By the general theory of equivariant Chow groups
(section \ref{chow}), we have approximations:
$$B\Or (V)= \stackrel{Lim}{m \rarr \infty} \ W_m/\Or (V),$$
$$B\SO (V)= \stackrel{Lim}{m \rarr \infty} \ W_m/\SO (V).$$
In this section, equivariant Chow rings of $\Or (k)$ and
$\SO (2k+1)$ are computed via the approximations
$$A^*_{\Or (V)}(\text{point})=  \stackrel{Lim}{m \rarr \infty} \ 
A^*(W_m/\Or (V)),$$
$$A^*_{\SO (V)}(\text{point})=  \stackrel{Lim}{m \rarr \infty} \ 
A^*(W_m/\SO (V)),$$
and the degeneracy loci results of section \ref{idealz}.

\subsection{The Chow Ring of $\Or (k)$}
Let $k\geq 1$.
Let  $V\eqq \com^k$ be equipped with a non-degenerate
quadratic form $Q$ preserved by $\Or (k)$. The quotient
$W_m/ \Or (k)$ can be explicitly realized as follows.
Let $\grass(k,m)$ be the Grassmannian of linear $k$-spaces
in $\com^m$.
Let $S \rarr \grass(k,m)$ be the tautological
sub-bundle. Let $Y_m \subset Sym^2 S^*$ be the
open locus of non-degenerate quadratic forms
on the fibers of $S$.
\begin{lm}
\label{fbb}
There is canonical $\Or (k)$-invariant
map $\tau: W_m \rarr Y_m$ which induces
an isomorphism $W_m/\Or (k) \eqq Y_m$.
\end{lm}
\bpf
Let $w\in W_m$. By the definitions,
$w$ naturally induces an injection $\iota_{w}:V \rarr \com^m$.
The quadratic form $Q$ then induces a non-degenerate
quadratic form $\iota_{w}(Q)$ on $\iota_{w}(V)$.
Let
$$\tau(w) = \iota_{w}(Q) \in Y_m.$$
It is easily checked that $\tau$ is an
algebraic morphism.
Let $g\in \Or (k)$. Then,
$$\iota_{g(w)}= \iota_{w} \circ g : V \rarr \com^m.$$
Hence, $\tau$ is $\Or (k)$-invariant.
Since the fibers of $\tau$ are exactly the
$\Or (k)$ orbits, the induced map
$$W_m/ \Or (k) \rarr Y_m$$
is a bijective morphism of nonsingular
complex algebraic varieties and thus an algebraic
isomorphism.
\epf

$Y_m$ is an approximation to $B\Or (k)$
up to codimension $m-k$.
$W_m\rarr Y_m$ is a principal $\Or(k)$-bundle
(see the Appendix).
The pull-back of the 
tautological sub-bundle $S \rarr \grass(k,m)$
to $Y_m$ is the vector bundle on $Y_m$
induced by the 
principal $\Or(k)$-bundle $W_m\rarr Y_m$
and the representation $V$.
The Chow ring of the Grassmannian $\grass(k,m)$
is freely generated by the Chern classes
$c_1, \ldots, c_k$ of $S$ up to codimension $m-k$
(the relations start in codimension $m-k+1$).
By Lemma \ref{paul}, the Chow ring of
$Y_m$ is isomorphic to
$$\Z[c_1, \ldots, c_k] \ / \ (\beta_1, \ldots, \beta_k)$$
up to codimension $m-k$ where
\begin{equation}
\label{bbbt}
\frac{c(S^*)}{c(S)}=
1+ \beta_1 +  \ldots + \beta_k+ \ldots.
\end{equation}
Induction and simple algebra establishes:
$$(\beta_1, \ldots, \beta_k)=(2c_1, 2c_3, 2c_5, \ldots).$$
The Chow ring limit $m\rarr \infty$ of $A^*(Y_m)$
is now easily seen to yield:
$$A^*_{\Or (k)}(\text{point})= 
\Z[c_1, \ldots, c_k] \ / \ (2c_1, 2c_3, 2c_5, \ldots).
$$
Theorem \ref{chor} is proven for $\Or (k)$.

\subsection{The Chow ring of $\SO (2k+1)$}
\label{ort}
Let $k\geq 0$.
Let $V\eqq \com^{2k+1}$ be equipped with a non-degenerate
quadratic form preserved by $$\SO (2k+1)\subset \GL (V).$$
Let $\com^* \subset \GL (V)$ be the scalars.
Since $V$ is odd dimensional $$\com^* \cap \SO (2k+1)= \{1\},
\ \ \com^* \times \SO(2k+1) \subset \GL(V).$$
The approximations $W_m$ to $E\SO (2k+1)$ are
used.
There is a natural free $\GL (V)$-action on $W_m$
which induces a free
 $\SO (2k+1)$-action and a free 
scalar $\com^*$-action on $W_m$.
The $\SO (2k+1)$-action and the $\com^*$-action
commute.
There is a commutative diagram:
\begin{equation}
\label{ffibb}
\begin{CD}
 W_m @>>>  W_m/ \SO(2k+1) \\
@VVV @VVV \\
 W_m/\com^* @>{\tau}>> W_m/\  \com^* \times \SO (2k+1)\\
\end{CD}
\end{equation}
All morphisms are group quotients:
the horizontal maps are free $\SO (2k+1)$-quotients,
the vertical maps are free $\com^*$-quotients.
See the Appendix for a discussion of these algebraic
quotient problems.

The quotients in diagram (\ref{ffibb}) are
analyzed. Let $S \rarr \grass(2k+1, m)$
be the tautological sub-bundle over the Grassmannian.
By an argument identical to Lemma \ref{fbb}, it
is seen that
$$W_m/ \ \com^* \times \SO (2k+1) \eqq Z_m$$
where $Z_m \subset \proj(Sym^2 S^*)$
is the locus of non-degenerate quadratic forms 
on the fibers of $S$.
Hence, $W_m/ \SO (2k+1) \rarr Z_m$
is a $\com^*$-bundle. Let $N \rarr Z_m$
be the line bundle associated to this $\com^*$-bundle.
On the left side of the diagram,
$$W_m/\com^* \subset \proj (\sumo_{1}^m V^*)$$
is the projective spanning locus.
$A^1( W_m/\com^*)= \Z$ and $W_m \rarr W_m/\com^*$
is the $\com^*$-bundle associated to the
generator $\oh_{\proj}(-1)$ of $A^1(W_m/\com^*)$.
\begin{lm}
\label{ddd}
$A^1(Z_m) \eqq \Z$ and $c_1(N)$ is a generator.
\end{lm}
\bpf
Consider the inclusion $Z_m \subset \proj(Sym^2 S^*)$.
Let 
$$\tau: W_m/\com^* \rarr Z_m \subset \proj(Sym^2 S^*)$$
be the natural map. Let $\barr{N}$ denote
an extension of $N$ to $A^1(\proj(Sym^2 S^*))$.
Since $\tau^*(\barr{N})= \oh_{\proj}(-1)$ 
generates $A^1(W_m/\com^*)\eqq \Z$,
the kernel $K$ of
$$ \tau^*: A^1(\proj(Sym^2 S^*)) \rarr A^1(W_m/ \com^*)$$
is isomorphic to $\Z$.
The class $[D]$ of the locus of degenerate
quadratic forms is in $K$.
$A^1(\proj(Sym^2 S^*))\eqq \Z c_1 \oplus \Z L$
where $c_1=c_1(S)$ and $L$ is 
the canonical class $\oh_{\proj}(1)$. The class of $[D]$ is
$-2c_1+ (2k+1) L$ which is not divisible in 
$A^1(\proj(Sym^2 S^*))$.
Hence $K$ is generated by $[D]$. Therefore
$\tau^*: A^1(Z_m) \rarr A^1(W_m/ \com^*)$ is
an isomorphism and $c_1(N)$ is a generator of $A^1(Z_m)$.
\epf 

There is now enough information to compute the
Chow ring of the approximation $W_m/\SO (2k+1)$ 
to $B\SO (2k+1)$.
As before, $W_m/\SO (2k+1)$ is an approximation up to codimension
$m-(2k+1)$.
The Chow ring of $\grass(2k+1, m)$
is freely generated by the Chern classes
$c_1, \ldots, c_{2k+1}$ of the tautological sub-bundle
$S$ up to codimension $m-(2k+1)$.
By Lemma \ref{pauly}, the
Chow ring of $Z_m$ (up to codimension $m-(2k+1)$) has
a presentation:
$$\Z[c_1, \ldots, c_{2k+1}, L]/ 
(p(L),\beta_1', \ldots, \beta_{2k+1}')$$
where $L$ is the class of $\oh_{\proj}(1)$, $p(L)$ is the
Chern polynomial satisfied by $L$, and 
$$\frac{c(S^*\otimes \oh_{\proj}(1))}{c(S)}=
1+ \beta'_1 + \ldots + \beta'_{2k+1}+ \ldots.$$
Finally, since $W_m/\SO (2k+1)$ is the
total space of the $\com^*$-bundle associated
to the line bundle $N\rarr Z_m$,
\begin{equation}
\label{ttoott}
\Z[c_1, \ldots, c_{2k+1}, L]/ 
(c_1(N),p(L),\beta_1', \ldots, \beta_{2k+1}')
\end{equation}
is a presentation of the Chow ring of $W_m/\SO (2k+1)$
(up to codimension $m-(2k+1)$).
Since $c_1(N)$ generates $A^1(Z_m)$ and the pair $\{c_1, L\}$
also generate $A^1(Z_m)$, (\ref{ttoott}) is 
equivalent to:
$$\Z[c_1, \ldots, c_{2k+1}, L]/ 
(c_1, L ,p(L),\beta_1', \ldots, \beta_{2k+1}').$$
By the defintions of
$p(L)$ and the elements $\beta'_i$,  there
is an equality of ideals 
$$(c_1,L, p(L), \beta'_1,\ldots, \beta_{2k+1}')
=
(c_1,L, c_s(Sym^2 S^*), \beta_1, \ldots, \beta_{2k+1})$$
where $s= rank(Sym^2 S^*)$ and the $\beta_i$ are determined
by (\ref{bbbt}).
\begin{lm}
$c_s(Sym^2 S^*) \in (\beta_1, \ldots, \beta_{2k+1})$.
\end{lm}
\bpf
Consider the total space $Sym^2 S^*$.
There is an isomorphism $A^*(Sym^2 S^*) \eqq 
A^*(\grass(2k+1,m))$.
The pull-back of the bundle $$Sym^2 S^* \rarr \grass(2k+1,m)$$
to the total space $Sym^2 S^*$ has a canonical section $\tau$.
The zero scheme of $\tau$ is contained in the
locus of degenerate quadratic forms $D\subset Sym^2 S^*$.
Also, the zero scheme of $\tau$ represents the class 
$$c_s(Sym^2 S^*)\in A^*(Sym^2 S^*).$$
Therefore,  $c_s(Sym^2 S^*) \in I_2$. The proof is
complete by Lemma \ref{paul}. \epf

\noindent As before,
$(\beta_1, \ldots, \beta_{2k+1})= (2c_1, 2c_3,2c_5, \ldots, 2c_{2k+1}).$
Hence, the Chow ring of $W_m/ \SO (2k+1)$ up to
codimension $m-(2k+1)$ has a presentation:
$$\Z[c_1, \ldots, c_{2k+1}]/ 
(c_1, 
2c_3, 2c_5,\ldots,  2c_{2k+1}).$$
The limit process yields Theorem \ref{chor} for $\SO (2k+1)$.

\section{\bf The Proof of Proposition \ref{alltor}}
\subsection{}
We follow the notation of section \ref{prezz}.
Let $V$ be a fixed $2$-dimensional $\com$-vector space.
Let $\proj^1\eqq \proj(V)$. 
Let $$U\subset \bigoplus_{0}^{d} Sym^d(V^*)$$
denote the non-degenerate locus parameterizing
bases of the linear series of $\oh_{\proj^1}(d)$
on $\proj^1$. $\GL(V)$ acts on $U$ properly with finite
stabilizers and geometric quotient (see the Appendix)
isomorphic to $H(d)$.
By Proposition \ref{dane},
\begin{equation}
\label{xeq}
A^*(d) \otimes_{\Z}{\Q} \eqq  A^*_{\GL (V)}(U) \otimes _{\Z}{\Q}.
\end{equation}
By the definition of equivariant Chow groups,
$$A^*_{\GL(V)}(U)= A^*(U \times^{\GL(V)} E \GL (V)).$$

\subsection{The Chow Rings of $\GL(V)$ and $\SL(V,n)$}
\label{slvn}
Algebraic approximations to $E\GL (V)$ are easily found.
Since related results about the groups $\SL(V,n)$ are
need in section \ref{ode}, a unified development
is presented here. Recall $SL(V,n) \subset \GL(V)$
is defined to be $det^{-1}(\Z/n\Z)$ where $det: \GL \rarr \com^*$
is the determinant homomorphism and 
$\Z/ n\Z$ is th group of $n^{th}$ roots of unity.

As in the orthogonal cases, the easiest approach to
$E \GL (V)$ and $E \SL (V,n)$ is via sums of the
representation $V^*$.
As before, let $m>>0$ and let 
$$W_m \subset \sumo_{1}^{m} V^*$$
be the spanning locus.
The induced $\GL (V)$ and $\SL (V,n)$-actions on
$W_m$ are free and have geometric quotients (see the
Appendix) which
approximate $B \GL (V)$ and $B \SL (V,n)$
up to codimension $m-2$.

It is easily seen that $W_m / \GL(V) \eqq \grass(2,m)$.
Since the Chow ring of this Grassmannian (up to
codimension $m-2$) is
freely generated by the Chern classes $c_1$ and $c_2$
of the tautological sub-bundle, 
$$A^*(W_m/ \GL (V)) \eqq \Z[c_1, c_2]$$
up to codimension $m-2$.
Taking the $m\rarr \infty $ limit, 
$$A^*_{\GL (V)}(\text{point}) \eqq \Z[c_1, c_2].$$
Similarly, $W_m/ \SL (V, n)$ is the total
space of the $n^{th}$ tensor power of the
line bundle $\bigwedge^2 S$ over $\grass(2,m)$.
Hence up to codimension $m-2$,
$$A^*(W_m/ \SL (V,n)) \eqq \Z[c_1, c_2]/ (nc_1).$$
Taking the $m \rarr \infty$ limit,
$$A^*_{\SL (V,n)}(\text{point}) \eqq \Z[c_1, c_2]/(nc_1).$$

\subsection{Proposition \ref{alltor}}
\label{abe}
The quotient $U \times ^{\GL (V)} E\GL (V)$
is analyzed via approximation. 
$V \times ^{\GL(V)} W_m$ is the tautological
sub-bundle $S$ over $\grass(2,m)$.
$$U \times ^{\GL(V)} W_m \subset \sumo_{0}^{d} Sym^d(V^*)
\times ^{\GL(V)} E$$ is the non-degenerate
open locus in the total space of the bundle
$\sumo_{0}^{d} Sym^d (S^*)$ over $\grass(2,m)$.
By Lemma \ref{petey}, there
is an isomorphism
$$A^*(U \times ^{\GL(V)} W_m) \eqq
\Z[c_1, c_2]/ (\alpha_1, \ldots, \alpha_{d+1})$$
up to codimension $m-2$ where
$$\frac{1}{c(Sym^d (S))}=
1+ \alpha_1+ \ldots + \alpha_{d+1} + \ldots.$$
The
ideal generated by $(\alpha_1, \ldots, \alpha_{d+1})$
is equal to the ideal generated by the 
first $d+1$ Chern classes of $Sym^d (S)$.

Taking the $m \rarr \infty$ limit, a presentation
of $A^*_{\GL (V)} (U)$ is obtained.
$A^*_{\GL (V)}(U)$ is generated (as a ring) 
in codimensions $1$, $2$ by elements $c_1$, $c_2$ respectively.
There are $d+1$ relations given as follows. Let $S$ be a rank
$2$ bundle with Chern classes $c_1$ and $c_2$. The $d+1$
Chern classes of $Sym^d(S)$ are the relations.
\begin{lm}
\label{ater}
$A^*_{\GL(V)}(U) \otimes \Q$ is zero is positive codimension.
\end{lm}
\bpf
A standard  calculation yields:
\begin{equation}
\label{see1}
c_1(Sym^d(S))=\frac{d(d+1)}{2} c_1,
\end{equation}
\begin{equation}
\label{see2}
c_2(Sym^d(S))= \frac{d(d-1)(d+1)(3d+2)}{24} c_1^2 +
\frac{d(d+1)(d+2)}{6}  c_2.
\end{equation}
Since the coefficients of $c_1$ and $c_2$ 
never vanish for positive $d$ in equations
(\ref{see1}) and (\ref{see2}) respectively,
the first two Chern classes of $Sym^d(S)$ generate
the ideal $(c_1, c_2)$ in $\Q[c_1,c_2]$.
\epf

\noindent
Lemma \ref{ater} and the isomorphism (\ref{xeq}) establish
Proposition \ref{alltor}.

\section{{\bf $A^*(d)$, $d$ Even}}
\label{evan}
The notation of section \ref{prezz} is used.
Let $d=2n$ (where $n\geq 1$).
Let $V\eqq \com^2$.
There is a free $\PGL (V)$-action
on $\proj(U) \subset \proj(\sumo_{0}^{d} Sym^d V^*)$
with geometric quotient (see the Appendix) 
isomorphic to $H(d)$.
By Proposition \ref{dane},
$$A^*(d) \eqq A^*_{\PGL (V)} (\proj (U)).$$
By the definition of equivariant Chow groups,
$$A^*_{\PGL (V)} (\proj (U)) \eqq A^*( \proj(U) \times
^{\PGL (V)} E \PGL (V)).$$
The Chow ring $A^*( \proj(U) \times
^{\PGL (V)} E \PGL (V))$ is computed in this section
for $d=2n$.

Consider the $3$-dimensional representation
$Sym^2(V)$ of $\PGL(V)$. This respresentation
leaves invariant a unique
(up to $\com^*$) quadratic form $Q$ on $Sym^2(V)$.
A group isomorphism
$\PGL (V) \eqq \SO (3)$ is induced by this quadratic
form.
The dual of the standard $3$-dimensional representation
of $\SO(3)$ corresponds to the representation
$Sym^2 (V^*)$ of $\PGL (V)$.
Let 
$$A_m \subset \sumo_{1}^{m} Sym^2 (V^*)$$
be the spanning locus. The approximations
$A_m/ \PGL (V)$ to $B \PGL (V)$ correspond
exactly to the approximations
$W_m/ \SO (3)$ to $B \SO (3)$ defined in section \ref{orthoo}.
$A_m/ \PGL(V)$ is therefore the total space of 
a $\com^*$-bundle $N \rarr Z_m$. 
$Z_m$ is the open set of non-degenerate
quadratic forms in
$\proj(Sym^2 (S^*))$ over the Grassmannian $\grass(3,m)$.
Let $B_m$ denote this approximation to $B \PGL (V)$.

$Sym^d(V^*)$ is a $\PGL(V)$ representation
for $d$ even ({\em not} for $d$ odd). Hence,
$$ Sym^d(V^*) \times^{\PGL(V)} A_m$$
is a rank $d+1$ vector bundle $F_d \rarr B_m$.
The quotient
$$\proj(U) \times
^{\PGL (V)} A_m \subset \proj(\sumo_{0}^{d} Sym^d(V^*))
\times ^{\PGL (V)} A_m$$
is simply the projective non-degenerate locus
in $\proj(\sumo_{0}^{d} F_d)$. 

The first step is to identify the bundle $F_d \rarr B_m$.
There is a tautological sub-bundle $S \rarr B_m$
obtained from the Grassmannian. There is a
tautological equivalence
$S^* \eqq F_2$. More generally, there
is a tautological sequence on $\proj(Sym^2 (S^*))$:
\begin{equation}
\label{toto}
 0 \rarr \oh_{\proj}(-1) \otimes Sym^{n-2} (S^*)
\rarr Sym^n (S^*) \rarr Q_n \rarr 0 
\end{equation}
for all $n\geq 2$.
Let $([q], P)\in \proj(Sym^2 (S^*))$ where
$P\subset \com^m$ is a linear $3$-space
and $0 \neq q \in Sym^2(P^*)$. The fiber
of $\oh_{\proj}(-1)$ over $([q],P)$ is
simply $\com \cdot q$.
The left inclusion in sequence (\ref{toto}) is
determined by the canonical multiplication map:
$$0 \rarr \com \cdot q \otimes Sym^{n-2} (P^*) \rarr 
     Sym^n(P^*).$$
Again, there is a tautological equivalence
$F_{2n} \eqq Q_n$ on $B_m$.

Note $A^1(B_m)=0$ by Lemma \ref{ddd}.
The Chern polynomial of $F_d$ on $B_m$ is therefore:
$$c(F_d)= \frac{c(Sym^n (S^*))}{c(Sym^{n-2}(S^*))}.$$
Now, by Lemma \ref{pete}, a presentation
of $A^*(\proj(U) \times^{\PGL(V)} B_m)$
up to codimension $m-3$ is obtained by
$$A^*(B_m)[\cal{L}]/
(p(\cal{L}), \alpha'_1, \ldots, \alpha'_{d+1})$$
where $\cal{L}$ is the class of 
$\oh_{\proj(\sumo_{0}^{d} F_d)}(1)$
and 
$$\frac{ (1+\cal{L})^{d+1}}{c(F^*_d)}= 
\frac{(1+\cal{L})^{d+1} \cdot c(Sym^{n-2} S)} {c(Sym^n S)} =
1+ \alpha'_1+ \ldots+ \alpha'_{d+1} 
\ldots .$$
By the presentation of $A^*(B_m)$ in section \ref{ort},
it follows
$$A^*(\proj(U) \times^{\PGL(V)} B_m) \eqq
 \Z[c_1, c_2, c_3, \cal{L}]/ 
(p(\cal{L}), c_1, 2c_3,  \alpha'_1, \ldots, \alpha'_{d+1})$$
up to codimension $m-3$. Taking the $m\rarr \infty$
limit, 
$$ A^*(d) \eqq A^*_{\PGL(V)}(\proj(U)) \eqq
\Z[c_1, c_2, c_3, \cal{L}]/ 
(p(\cal{L}), c_1, 2c_3,  \alpha'_1, \ldots, \alpha'_{d+1}).$$
The relation $p(\cal{L})$ is of codimension $(d+1)^2$.
Since the dimension of $H(d)$ is $d^2+ 2d-3=(d+1)^2-4$ and
the generators $c_1,c_2, c_3$, and $\cal{L}$ have
dimension at most 3, $p(\cal{L})$
is a relation among classes that are already zero.
Hence,
\begin{equation}
\label{freedy}
A^*(d) \eqq A^*_{\PGL(V)}(\proj(U)) \eqq
\Z[c_1, c_2, c_3, \cal{L}]/ 
(c_1, 2c_3,  \alpha'_1, \ldots, \alpha'_{d+1}).
\end{equation}

Following [EG], there is natural
map 
$$\proj(U)\times^{\PGL(V)} A_m \rarr \proj(U)/ \PGL(V) \eqq H(d)$$
which expresses $\proj(U)\times^{\PGL(V)} A_m$ as an open set
of a vector bundle over $H(d)$. This fibration
induces an isomorphism on Chow rings (up to
codimension $m-3$). The classes $c_i\in A^*(d)$ are
easily identified via this isomorphism (up to
codimension $m-3$):
\begin{equation}
\label{ddod}
A^*(\proj(U)\times^{\PGL(V)} A_m) \eqq A^*(d).
\end{equation}
They are the Chern classes of the vector
bundle obtained from the
the principal $\PGL(V)$-bundle $\proj(U) \rarr
H(d)$ and the representation $Sym^2(V)$.
Let $\cal{H}\in A^1(d)$ be the class of 
curves meeting a fixed codimension 2 linear
space $P$ of $\proj^d$. $\cal{H}$ corresponds via the isomorphism
(\ref{ddod}) to a resultant
class in $A^1(\proj(U)\times^{\PGL(V)} A_m)$.
Routine calculations show $\cal{H}=2d \cal{L}$
where $2d$ is the degree of the resultant of
degree $d$ polynomials.
Since $(d+1)\cal{L}=0$ by the presentation (\ref{freedy}).
$$n \cal{H} = d^2 \cal{L} = (d^2-1+1)\cal{L} = (d-1)(d+1)\cal{L} 
+\cal{L} = \cal{L}.$$
The proof of Theorem \ref{evan} is complete.

\section{{\bf $A^*(d)$, $d$ Odd}}
\label{ode}
Let $d=2n-1$ (where $n\geq 1$).
Let $V\eqq \com^2$. 
There is a canonical, $\GL (V)$-equivariant, multilinear map
$$\mu:\bigoplus_0^{2n-1} 
Sym^{2n-1}(V^*) \rarr \bigwedge^{2n} Sym^{2n-1}(V^*)$$
given by the exterior product (see section \ref{prezz}).
\begin{lm} The
$\SL (V,n)$-action on 
$\bigwedge^{2n} Sym^{2n-1}(V^*)$ is trivial. 
\end{lm}
\bpf
Since the  $1$-dimensional representations of
$\SL(V)$ are trivial, the action of $\SL(V)$
on
$\bigwedge^{2n} Sym^{2n-1}(V^*)$ is certainly trivial. Let
$H\subset\SL(V,n)$ be the subgroup of scalars.
$H$ is the multiplicative group of scalar $2n^{th}$
roots of unity, 
Let $\xi \in H$ be a scalar. $\xi$ acts
on $\bigwedge^{2n} Sym^{2n-1}(V^*)$ by the
scalar $\xi^{(2n)(2n-1)}= 1$.
It is easily checked that $\SL(V,n)$
is generated (as a group) by $H$ and $\SL(V)$.
Hence, the $\SL(V,n)$-action is trivial.
\epf

\noindent
Let 
$Y= \mu^{-1} (p)$ where $0 \neq p \in\bigwedge^{2n} Sym^{2n-1}(V^*)$.
There is an $\SL (V,n)$-action on $Y$. 
\begin{lm}
\label{freeaq}
The $\SL (V,n)$-action on $Y$ is free with geometric
quotient $H(d)$.
\end{lm}
\bpf
Certainly $\SL(V,n)$ acts on $Y$ since the
$\SL(V,n)$-action on $\bigwedge^{2n} Sym^{2n-1}(V^*)$
is trivial.
Let 
$$U \subset \sumo_0^{2n-1} Sym^{2n-1} (V^*)$$
be the non-degenerate locus.
First, it is shown that the $\SL(V,n)$-action
on $U$ is free. Since $Y\subset U$, $\SL(V,n)$
acts freely on $Y$.

Let $u\in U$. Suppose $g \in \SL(V,n)$ satisfies
$g\cdot u = u$. $\PGL(V)$ acts freely
on $\proj(U)$. Let $\pi: \SL(V,n) \rarr \PGL(V)$.
Then,
$$\pi(g) \cdot \proj(u)= \proj(u).$$
Hence, $\pi(g) = 1 \in \PGL (V)$. The
element $g$ is therefore a scalar in $\SL(V,n)$
equal to a $2n^{th}$ root of unity $\xi$.
Then, $g$ acts on $u$ by
the scalar $\xi^{2n-1}$. Since
$g \cdot u=u$, $\xi^{2n-1} =1$.
Since $(2n, 2n-1)=1$, $\xi^{2n-1}=1$ implies
$\xi=1$. Therefore, $g=1 \in \SL(V,n)$.
The $\SL(V,n)$-action on $U$ is free.

It is now shown the quotient $Y/ \SL(V,n)$ is isomorphic
to $H(d)$. There are natural, equivariant, algebraic projection
maps:
$$Y \rarr \proj(U),$$
$$\pi: \SL(V,n) \rarr \PGL(V).$$
These maps induce a natural surjective map on quotients:
$$\phi: Y/ \SL(V,n) \rarr \proj(U)/ \PGL(V) \eqq H(d).$$
It suffice to prove $\phi$ is injective. (A bijective
map of nonsingular complex algebraic varieties
is an algebraic isomorphism.)

Let $y_1, y_2 \in Y$ be points. Let
$[y_1], [y_2]\in \proj(U)$ denote the corresponding
points. Suppose there exists an element $\gamma \in \PGL(V)$
satisfying $\gamma \cdot [y_1]=[y_2]$.
To prove $\phi$ is injective,
it must be shown that  $y_1$ and $y_2$ are in the same
$\SL(V,n)$ orbit.

Let $g\in \SL(V,n)$ satisfy $\pi(g)=\gamma$.
Then,
 $[g \cdot y_1] =[y_2]$.
Hence $g \cdot y_1= (\lambda y_2)$
where $\lambda \in \com^*$ is a scalar.
By the conditions $g \cdot y_1, y_2 \in Y$,
it follows $\lambda^{2n}=1$.
Since $(2n, 2n-1)=1$, a $2n^{th}$ root
of unity $\xi \in \SL(V,n)$ can be found satisfying
$\xi^{2n-1}=\lambda^{-1}$.
Let $h \in \SL(V,n)$ be determined by $h= \xi \cdot g$.
$$h\cdot y_1= \xi \cdot g \cdot y_1= \xi \cdot (\lambda y_2)
= \xi^{2n-1}(\lambda y_2)= y_2.$$
Therefore $y_1$ and $y_2$ are in the same $\SL(V,n)$ orbit.
\epf

\noindent
There is a canonical
isomorphism of graded rings
$$A^*(d=2n-1)  \eqq  A^*_{\SL (V,n)}(Y).$$ 
The equivariant Chow ring $A^*_{\SL (V,n)}(Y)$
is computed in this section.

The approximations $W_m$ and $W_m/ \SL(V,n)$ to
$E \SL(V,n)$ and $B \SL(V,n)$ determined in
section \ref{slvn} are used here.
Recall $W_m/ \SL(V,n) \rarr \grass(2,m)$
is the $\com^*$-bundle associated to the $n^{th}$ 
tensor power of $\wedge^2 S$ (where $S$ is
the tautological sub-bundle over $\grass(2,m)$).
Since $$Y \subset U \subset \sumo_{0}^{2n-1} Sym^{2n-1} (V^*),$$
there are inclusions:
$$ Y \times^{\SL(V,n)} W_m \ \subset \
   U \times ^{\SL(V,n)} W_m \ \subset \
   \sumo_{0}^{2n-1} Sym^{2n-1}(V^*) \times ^{\SL(V,n)} W_m.$$
Let $F_n=V^* \times^{\SL(V,n)} W_m$. $F_n$ is an algebraic
vector bundle over $W_m/ \SL(V,n)$. $F_n$ is
easily identified as the pull-back of $S^*$ to
$W_m/ \SL(V,n)$.
$U \times ^{\SL(V,n)} W_m$ is the
affine non-degenerate locus (i) of section \ref{idealz}
associated to the bundle $Sym^{2n-1} F_n\eqq Sym^{2n-1} S^*$.

\begin{lm}
\label{qwq}
There is an isomorphism 
$$ \epsilon: \com^* \times (Y \times^{\SL(V,n)} W_m) 
\eqq U \times ^{\SL(V,n)} W_m.$$
\end{lm}
\bpf
Let $\SL(V,n)$ act trivially on $\com^*$.
Define a $\SL(V,n)$ equivariant
isomorphism
$$\delta: \com^* \times Y \rarr U$$
by the following:
$$\delta\big( \lambda, (\omega_0, \omega_1\ldots, \omega_{2n-1})\big)
= (\lambda \omega_0, \omega_1, \ldots, \omega_{2n-1}).$$
The isomorphism $\delta$ induces
isomorphisms:
$$\com^* \times (Y \times^{\SL(V,n)} W_m) \eqq
(\com^* \times Y) \times^{\SL(V,n)} W_m \eqq
 U \times^{\SL(V,n)} W_m.$$
Let $\epsilon$ be the composition.
\epf

$W_m/ \SL(V,n)$ approximates $B \SL(V,n)$
up to codimension $m-2$.
The Chow ring of $Y \times ^{\SL(V,n)} W_m$
is now computed (up to codimension $m-2$).
By Lemma \ref{qwq}, there is an isomorphism:
$$A^*(Y \times ^{\SL(V,n)} W_m) \eqq 
A^*(U \times ^{\SL(V,n)} W_m).$$
Since $U$ is the affine non-degenerate
locus associated to the bundle 
$$Sym^{2n-1} S^* \rarr (W_m / \SL(V,n)),$$
Lemma \ref{pete} can be applied.
Recall the Chow ring of $W_m / \SL(V,n)$ (up to
codimension $m-2$)
has a presentation $\Z[c_1, c_2]/ (nc_1)$.
Hence, there is an isomorphism (up to
codimension $m-2$):
$$A^*(Y \times ^{\SL(V,n)} W_m) \eqq
\Z[c_1, c_2]/ (nc_1, \alpha_1,  \ldots, \alpha_{d+1})$$
where
$$\frac{1}{c(Sym^d S)}=
1+ \alpha_1 + \ldots + \alpha_{d+1}+ \ldots.$$
The proof of Theorem \ref{hodd} is complete.

\section {\bf Examples}
Since $H(1)$ is a point, $A^*(1)$ is the trivial
$\Z$-algebra (which agrees with the presentation
of Theorem 2).
$H(2)$ is the space of nonsingular plane conics.
By Theorem 1, $A^*(2)$ is generated
by $c_2$, $c_3$, and $\cal{L}=\cal{H}$ subject to 
$4$ relations:
$$2c_3=0,$$
$$3\cal{H}=0, \ -c_2+3\cal{H}^2=0, \ -c_3+\cal{H}^3=0.$$
Since $\cal{H}$ is 3-torsion, $c_2=0$. 
Since $c_3$ is two torsion, the last
equation can be reduced to $c_3=\cal{H}^3=0$.
Therefore $A^*(2)$ is given by
$$\Z[\cal{H}]/ (3\cal{H}, \cal{H}^3).$$
Since $H(2)$ is an open set of the projective
space of plane conics, another approach to $A^*(2)$
is possible.
The class $\cal{H}$ is simply the restriction
of the hyperplane class which necessarily generates $H(2)$.
The relation $3\cal{H}$ 
can be obtained from the degree $3$ degeneracy locus
of singular plane conics.
The relation $\cal{H}^3$ is a 
consequence of the fact that the locus of
conics singular at a {\em fixed} point in $\proj^2$
is a linear $\proj^2$ in the $\proj^5$ of conics.

Let $d=3$, $n=2$, $d=2n-1$. By Theorem 2,
$A^*(3)$ is generated by $c_1$ and $c_2$ 
with relations:
$$ 2c_1=0,$$
$$6c_1=0, \ 11c_1^2+ 10c_2=0, \ 6c_1^3+30c_1c_2=0, 
\ 18c_1^2c_2+9c_2^2=0.$$
These relations simplify to yield the presentation:
$$A^*(3)= \Z[c_1,c_2]/ (2c_1, c_1^2+10c_2, c_1^3, c_1^2c_2,c_2^2).$$
In particular, $A^i(3)=0$ for $i\geq 4$.

\section {\bf Appendix On Algebraic Quotients}
\label{appx}
Let $\com$ be the ground field of complex numbers.
The geometric invariant theory terminology of [MFK]
is used here.
Let $\G$ be a reductive linear algebraic group.
A group action $\G \times X \rarr X$ is {\em proper}
if the natural map
$$\Psi: \G \times X \rarr X \times X$$
(given by the action and projection onto the second
factor) is a proper morphism.
The main result needed is the following:

\begin{pr}
Let $X$ be a quasi-projective variety with
a linearized $\G$-action satisfying $X^{stable}_{(0)}
=X$.
\label{qquot}
Then, the
 $\G$-action on $X$ is proper and 
there is a quasi-projective 
geometric quotient $X \rarr X/\G$.
\end{pr}
\bpf
Properness of the action is exactly Corollary 2.5
of [MFK]. The geometric quotient
is the main construction in geometric invariant
theory (Theorem 1.10 of [MFK]).
\epf

\noindent
The stable locus $X_{(0)}^{stable}$ is detected
by the Numerical Criterion.

Let $V\eqq \com^k$ be a vector space equipped
with a quadratic form.
All of the linear algebraic groups
considered in this paper are reductive:
$\GL(V)$, $\SL(V)$, $\PGL(V)$, $\SL(V,n)$, $\Or(V)$, $\SO (V)$,
$\com^*$, $\com^* \times \SO(V)$.
Let 
$$Span_m(V,d) \subset \sumo _{1}^{m} Sym^d(V*)$$
be the spanning locus (the locus of $m$-tuples
of vector of $Sym^d(V^*)$ which span $Sym^d(V^*)$).
The spanning loci
$U \subset \sumo_{0}^{d} Sym^d(\com^{2*})$
and
$W_m \subset \sumo_{1}^{m} V^*$
are special cases of $Span_m(V,d)$. 
The group actions considered in the paper are
of three forms:
\begin{enumerate}
\item[(i)] The natural $\G$-action on $X=Span_m(V,d)$
           where $\G \subset \GL(V)$ is
           a reductive subgroup. 
\item[(ii)] The $\G$-action on a $\G$-invariant
            subvariety $Y\subset Span_m(V,d)$ where
            $\G \subset \GL(V)$ is a reductive 
            subgroup.
\item[(iii)] The natural $\PGL(V)$-action on $X=\proj(Span_m(V,d))$.
\end{enumerate}
For example, the $SL(V,n)$-action
on $Y$ considered in section \ref{ode} is of form (ii).
Consider first (i) and (ii).
A linearization of the $\GL(V)$-action
can be found on $X$ satisfying  $X^{stable}_{(0)}=X$.
Such a linearization is found in section 1 of 
[P]. Since the stable locus
is detected by the Numerical Criterion, the
result for $\GL(V)$ implies $X^{stable}_{(0)}=X$
for the induced action of any reductive
subgroup $\G \subset \GL(V)$.
It is similarly simple to find a linearization
in case (iii) satisfying $X^{stable}_{(0)}=X$.
Therefore, Proposition \ref{qquot} applies
to all the quotient problems in the paper.

In [MFK], the $\G$-action $\G \times X \rarr X$
is defined to be {\em free} if
the natural map $\Psi$
is a closed embedding.
An action is {\em set-theoretically free}
if the stabilizers are trivial.
For the set-theoretically free 
actions considered in this paper,
the following Lemma is utilized.
\begin{lm}
Let $X$ be nonsingular.
Let $\G\times X \rarr X$ be a proper action.
In this case, set-theoretically free implies free.
\end{lm}
\bpf
Let $I \subset X\times X$ be the image of $\Psi$.
$I$ is a closed subvariety since $\Psi$ is proper.
It must be shown that $\G \times X \rarr I$
is an isomorphism. 
First it is shown that $I$ is nonsingular.
For this, it suffices to prove the
differential of $\Psi$ is injective at each point
of $\G \times X$. It is well known (over $\com$) that
set-theoretically trivial stabilizers are
also scheme-theoretically trivial. From this, 
the injectivity of the differential $d\Psi$ is easily deduced.
Now $\Psi: \G \times X \rarr I$ is
a bijective map of nonsingular complex
algebraic varieties and thus an algebraic
isomorphism.
\epf

\noindent A result of [MFK] (Proposition 0.9)
relates free quotients to principal $\G$-bundles.
\begin{pr}
Let $\G \times X \rarr X$ be an algebraic
group action with geometric quotient $X \rarr Y$.
If the action is free, then $X \rarr Y$ is
a (\'etale locally trivial) principal $\G$-bundle.
\end{pr}
\noindent
A principal $\G$-bundle $X \rarr Y$  and a representation
$\G\rarr\GL$ together yield a principal $\GL$-bundle
$X \times^{\G} \GL \rarr Y$.
Since every principal $\GL$-bundle is 
Zariski locally trivial ($\GL$ is a {\em special}
group in the sense of Grothendieck (see
{\em Anneau de Chow et Applications}, Seminaire
Chevalley 1958)), an algebraic vector bundle
over $Y$ is obtained.

\vspace{+10 pt}
\noindent
Department of Mathematics \\
University of Chicago \\
5734 S. University Ave. 60637 \\
rahul@math.uchicago.edu

\begin{thebibliography}{[MFK]}
\bibitem[EG]{a} D. Edidin and W. Graham, {\em Equivariant Intersection
                Theory}, preprint 1996.
\bibitem[EG2]{aa} D. Edidin and W. Graham, {\em Characteristic classes
of principal bundles in algebraic geometry}, preprint 1995.
\bibitem[F]{b} W. Fulton, {\em Intersection Theory}, 
             Springer-Verlag:
             Berlin, 1984.
\bibitem[MFK]{c} D. Mumford, J. Fogarty, and F. Kirwan, {\em
                 Geometric Invariant Theory}, Springer-Verlag: Berlin,
                 1994.
\bibitem[P]{d} R. Pandharipande, {\em The Chow Ring of the
 Non-linear Grassmannian}, preprint 1996.
\bibitem[Pr]{e} P. Pragacz, {\em Enumerative Geometry
of Degeneracy Loci}, Ann. Scient. Ec. Norm. Sup. 
{\bf 21}, 1988, p.413-454.
\bibitem[T] {f} B. Totaro, {\em The Chow Ring of the
Symmetric Group}, preprint 1994.
\end{thebibliography}
\end{document}